\documentclass[12pt]{amsart}
\usepackage{amsmath,amscd,amssymb,amsfonts}
\newcommand{\h}{\hbox}
\newcommand{\q}{\quad}
\newcommand{\nin}{\par\noindent}
\newcommand{\bs}{\par\bigskip}
\newcommand{\ms}{\par\medskip}
\newcommand{\sk}{\par\smallskip}
\newcommand{\bl}{\bigl}
\newcommand{\br}{\bigl}
\newcommand{\ssb}{\raise.15ex\h{${\scriptscriptstyle\bullet}$}}
\newcommand{\scc}{\,\raise.15ex\h{${\scriptstyle\circ}$}\,}
\newcommand{\msum}{\h{$\sum$}}
\newcommand{\mcap}{\h{$\bigcap$}}
\newcommand{\mcup}{\h{$\bigcup$}}
\newcommand{\mopl}{\h{$\bigoplus$}}
\newcommand{\mprod}{\h{$\prod$}}
\newcommand{\mcoprod}{\h{$\coprod$}}
\newcommand{\mtim}{\h{$\times$}}
\newcommand{\al}{\alpha}
\newcommand{\be}{\beta}
\newcommand{\C}{{\mathbf C}}
\newcommand{\D}{{\mathbf D}}
\newcommand{\Dc}{{\mathcal D}}
\newcommand{\E}{\widetilde{E}}
\newcommand{\e}{\varepsilon}
\newcommand{\F}{{\mathcal F}}
\newcommand{\f}{\tilde{f}}
\newcommand{\Ga}{\Gamma}
\newcommand{\g}{\widetilde{\gamma}}
\newcommand{\Hc}{{\mathcal H}}
\newcommand{\HH}{{\mathbf H}}
\newcommand{\I}{{\mathcal I}}
\newcommand{\iti}{\tilde{i}}
\newcommand{\jti}{\,\widetilde{\!j}}
\newcommand{\LL}{{\mathbf L}}
\newcommand{\La}{\Lambda}
\newcommand{\M}{{\mathcal M}}
\newcommand{\MM}{{\mathcal M}^{\ssb}}
\newcommand{\N}{{\mathbf N}}
\newcommand{\Oc}{{\mathcal O}}
\newcommand{\PP}{{\mathbf P}}
\newcommand{\pit}{\widetilde{\pi}}
\newcommand{\Q}{{\mathbf Q}}
\newcommand{\Qt}{\widetilde{Q}}
\newcommand{\R}{{\mathbf R}}
\newcommand{\So}{\,\overline{\!S}}
\newcommand{\Sc}{{\mathcal S}}
\newcommand{\Si}{\Sigma}
\newcommand{\si}{\sigma}
\newcommand{\sit}{\widetilde{\si}}
\newcommand{\spe}{{\rm sp}}
\newcommand{\sss}{{\mathbf s}}
\newcommand{\U}{{\mathcal U}}
\newcommand{\X}{\widetilde{X}}
\newcommand{\Xc}{{\mathcal X}}
\newcommand{\Y}{{\mathcal Y}}
\newcommand{\Yt}{\widetilde{Y}}
\newcommand{\Yct}{\widetilde{\mathcal Y}}
\newcommand{\Z}{{\mathbf Z}}
\newcommand{\Zc}{{\mathcal Z}}
\newcommand{\BM}{{\rm BM}}
\newcommand{\CH}{{\rm CH}}
\newcommand{\DR}{{\rm DR}}
\newcommand{\Gr}{{\rm Gr}}
\newcommand{\Hom}{{\rm Hom}}
\newcommand{\MHM}{{\rm MHM}}
\newcommand{\MHS}{{\rm MHS}}
\newcommand{\Sing}{{\rm Sing}}
\newcommand{\vir}{{\rm vir}}
\newcommand{\pc}{\prime\circ}
\newcommand{\into}{\hookrightarrow}
\newcommand{\too}{\longrightarrow}
\newcommand{\simto}{\buildrel\sim\over\longrightarrow}
\begin{document}
\title{Hirzebruch-Milnor Classes of Complete Intersections}
\author[L. Maxim ]{Laurentiu Maxim}
\address{L. Maxim : Department of Mathematics, University of
Wisconsin-Madison, 480 Lincoln Drive, Madison WI 53706-1388 USA}
\email{maxim@math.wisc.edu}
\author[M. Saito ]{Morihiko Saito}
\address{M. Saito: RIMS Kyoto University, Kyoto 606-8502 Japan}
\email{msaito@kurims.kyoto-u.ac.jp}
\author[J. Sch\"urmann ]{J\"org Sch\"urmann}
\address{J.  Sch\"urmann : Mathematische Institut, Universit\"at
M\"unster, Einsteinstr. 62, 48149 M\"unster, Germany}
\email{jschuerm@uni-muenster.de}
\dedicatory{To the memory of Friedrich Hirzebruch}
\begin{abstract}
We prove a new formula for the Hirzebruch-Milnor classes of global
complete intersections with arbitrary singularities describing the
difference between the Hirzebruch classes and the virtual ones.
This generalizes a formula for the Chern-Milnor classes in the
hypersurface case that was conjectured by S.~Yokura and was proved by
A.~Parusinski and P.~Pragacz. It also generalizes a formula of J.~Seade
and T.~Suwa for the Chern-Milnor classes of complete intersections with
isolated singularities.
\end{abstract}
\maketitle

\centerline{\bf Introduction}
\bs\nin
For the proof of his Riemann-Roch theorem [Hi], F.~Hirzebruch defined
the $\chi_y$-genus of a compact complex manifold $X$ by
$$\chi_y(X):=\msum_p\,\chi(\Omega_X^p)\,y^p\in\Z[y].$$
This specializes respectively to the Euler characteristic, the
arithmetic genus, and the signature of $X$ at $y=-1,0,1$.
It is the highest degree part of the cohomology Hirzebruch
characteristic class $T_y^*(TX)$ of the tangent bundle $TX$.
Using the Chern roots $\{\al_i\}$ for $TX$, this cohomology Hirzebruch
class is defined by 
$$\aligned T^*_y(TX)&:=\mprod_{i=1}^{\dim X}Q_y(\al_i)\in
\HH^{\ssb}(X)[y],\raise-7pt\h{}\\
\h{with}\q Q_y(\al)&:={\al(1+y)/\bl(1-e^{-\al(1+y)}
\br)}-\al y\in\Q[y][[\al]],\endaligned$$
where $Q_y(\al)$ is as in [Hi], 1.8 (see also (1.1) below), and
$\HH^k(X)=H^{2k}(X,\Q)$ in this paper.
(It is known that the theory works also for the Chow cohomology groups
as defined in [Fu1], [Fu2].)
Substituting $y=-1$, $0$, $1$, we see that
$Q_y(\al)\in\Q[y][[\al]]$ becomes respectively
$$1+\al,\q\al/(1 -e^{-\al}),\q\al/\tanh\al,$$
and hence $T^*_y(TX)$ specializes respectively to the Chern class
$c^*(TX)$, the Todd class $td^*(TX)$, and the Thom-Hirzebruch
$L$-class $L^*(TX)$, see [HBJ], Sect.~5.4.
In the smooth case, this cohomology class $T_y^*(TX)$ is identified
by Poincar\'e duality with the (Borel-Moore) homology class
$T_y^*(TX)\cap[X]$ which will be denoted by $T_{y*}(X)$.
\sk
A generalization of $T_{y*}(X)$ to the singular case was given by
[BSY] using the Du Bois complex in [DB] or $\Q_{h,X}\in D^b\MHM(X)$,
the bounded complex of mixed Hodge modules on $X$ whose underlying
$\Q$-complex is the constant sheaf $\Q_X$, where $\MHM(X)$ is the
category of mixed Hodge modules [Sa2].
The Hirzebruch class $T_{y*}(X)$ is actually a special case of
$T_{y*}(\M^{\ssb})$ defined for any $\M^{\ssb}\in D^b\MHM(X)$ in [BSY],
see (1.2) below.
\sk
Hirzebruch [Hi] also introduced the notion of virtual $\chi_y$-genus
(or $\chi_y$-characteristic) which is the $\chi_y$-genus of smooth
complete intersections $X$ in smooth projective varieties $Y$.
The {\it virtual Hirzebruch characteristic class} $T^{\,\vir}_{y*}(X)$
can be defined like the virtual genus (even in the singular case),
see (1.3) below.
In this paper, we show that the difference between the Hirzebruch class
and the virtual one is given by the {\it Hirzebruch-Milnor class}
supported on the singular locus of $X$, and prove an {\it inductive}
formula to {\it calculate} it explicitly in the case of global
complete intersections with arbitrary singularities as follows.
\sk
Let $Y^{(0)}$ be a smooth projective variety, and $L$ be a very ample
line bundle on $Y^{(0)}$. We have the graded ring
$$R_L:=\mopl_{k\ge 0}\,R_{L,k}\q\h{with}\q
R_{L,k}:=\Ga(Y^{(0)},L^{\otimes k}).$$
Let $s_j\in R_{L,a_j}$ with $\{a_j\}$ a decreasing
sequence of positive integers.
Let $r$ be a positive integer strictly smaller than $\dim Y^{(0)}$.
For $j\in[1,r]$, set
$$Y^{(j)}:=\mcap_{k\in[1,j]}\,s_k^{-1}(0)\subset Y^{(0)},$$
and assume ${\rm codim}_{Y^{(0)}}Y^{(j)}=j$. Define
$$X:=Y^{(r)},\q Y:=Y^{(r-1)}.$$
Note that $Y$ is not necessarily smooth. It is a complete intersection
of codimension $r-1$ in a smooth projective variety $Y^{(0)}$.
The sections $s_j\,(j\in[1,r])$ generate a graded ideal
$$I_X=\mopl_{k\ge 0}\,I_{X,k}\subset R_L.$$
Set
$$\PP(R_{L,k}):=(R_{L,k}\setminus\{0\})/\C^*,\q
\PP(I_{X,k}):=(I_{X,k}\setminus\{0\})/\C^*.$$
Replacing the $s_j$ without changing $I_X$, we may assume by
Proposition~(3.2) below
$$\Sing\,Y^{(j)}\subset Y^{(r)}=X\q(j\in[1,r]).$$
This condition is satisfied (and the $s_j$ generate the ideal $I_X$)
if the $s_j$ are sufficiently general, i.e.\ if $([s_j])$ belongs to a sufficiently small Zariski-open subset
$$\U_I\subset\mprod_{j=1}^r\,\PP(I_{X,a_j}),$$
where $[s_j]$ denotes the image of $s_j$ in $\PP(I_{X,a_j})$.
\sk
Take a sufficiently general $s'_r\in R_{L,a_r}$, and set
$$X':=s_r^{\prime -1}(0)\cap Y.$$
This is viewed as a complete intersection of codimension $r-1$ in the
{\it smooth} hypersurface section $s_r^{\prime -1}(0)$ of $Y^{(0)}$.
Here ``sufficiently general'' means that $[s'_r]$ belongs to a
sufficiently small Zariski-open subset $U_r$ of $\PP(R_{L,a_r})$.
It satisfies at least the following two conditions:
\sk\nin
(1) The hypersurface section $s_r^{\prime -1}(0)$ transversally
intersects all the strata of an algebraic Whitney stratification of $Y$
such that $Y\setminus X$ is a stratum.
\sk\nin
(2) For an embedded resolution of $X\subset Y$, the pullback of $X'$
is smooth and transversally intersects any intersections of the
irreducible components of the pullback of $X$.
\sk
Conditions~(1) and (2) respectively define a non-empty Zariski-open
subset of $\PP(R_{L,a_r})$ for each choice of an algebraic Whitney stratification of $Y$ or an embedded resolution of $(Y,X)$.
Note that the intersection of finitely many non-empty Zariski-open subsets
of $R_{L,a_r}$ is non-empty, and this is useful when we have to show the
independence of the choices of the $s_j$, $s'_j$ later.
Note also that the Thom $a_f$-condition follows from condition~(1)
(see [BMM], [Pa]), and condition~(2) follows from (1) if there
is an algebraic Whitney stratification of the embedded resolution such
that the desingularization is a stratified morphism, i.e.
any stratum of the resolution is smooth over a stratum of $Y$.

We have
$$\Si:=\Sing\,X\supset\Sing\,Y\supset\Si':=\Sing\,X'=\Sing\,Y\cap X'.$$
Let $i_{A,B}$ denote the inclusion morphism for $A\subset B$ in
general.
Set $f:=(s_r/s'_r)|_{Y\setminus X'}$.
\ms\nin
{\bf Theorem~1.} {\it With the above notation and assumptions,
there is $M_y(X)\in\HH_{\ssb}(\Si)[y]$, called the Hirzebruch-Milnor
class, where $\HH_{\ssb}(\Si):=H_{2\ssb}^{\BM}(\Si,\Q)$, and satisfying
$$T^{\,\vir}_{y*}(X)-T_{y*}(X)=(i_{\Si,X})_*M_y(X),
\leqno(0.1)$$
\vskip-20pt
$$M_y(X)=T_{y*}\bl((i_{\Si\setminus X',\Si})_{!\,}\varphi_f
\Q_{h,Y}\br)+(i_{\Si',\Si})_*M_y(X').
\leqno(0.2)$$
More precisely, choosing sufficiently general
$s'_j\in R_{L,a_j}\,(j\in[1,r])$ inductively,
we have a mixed Hodge module $\M(s'_1,\dots,s'_r)$ on $\Si$ such that
$$M_y(X)=(-1)^{\dim X}T_{y*}\bl(\M(s'_1,\dots,s'_r)\br),$$
and $M_y(X)$ is independent of the choice of the $s'_j$.}
\ms
Here $\varphi_f\Q_{h,Y}$ is viewed as a mixed Hodge module
on $\Si\setminus X'$ up to a shift.
(Note that $\Q_Y[\dim Y]$ is a perverse sheaf for any local complete
intersection $Y$.)
The monodromy around $X'$ of the local systems $H^j\varphi_f\Q_{h,Y}|_S$
on a stratum $S$ of a Whitney stratification coincides with the
Milnor monodromy, since $f=(s_r/s'_r)|_{Y\setminus X'}$ by definition.
So it is rare that the monodromy is trivial, and hence the situation is
quite different from the one in [CMSS].
By Theorem~1, we have $T^{\,\vir}_{y,k}(X)=T_{y,k}(X)$ for $k>\dim\Si$
where $T^{\,\vir}_{y,k}(X)\in\HH_k(X)[y]$ is the degree $k$ part of
$T^{\,\vir}_{y*}(X)$.
This generalizes [CMSS], Cor.~3.3 in the hypersurface case.
Theorem~1 was proved in [CMSS], (1.17) in the case of hypersurfaces with
isolated singularities.
\sk
The above $M_y(X)$ may depend on the choice of the $s_j$.
However, it is well-defined if $([s_j])$ belongs to a sufficiently
small non-empty Zariski-open subset $\U_I$ of
$\prod_{j=1}^r\PP(I_{X,a_j})$ (although it is not necessarily easy to
write down explicitly $\U_I$), see Proposition~(4.4) below.
Note that $M_y(X)$ coincides with this canonical one if it is
invariant by any small perturbation of the $s_j$, i.e. if it is constant
for any element in a sufficiently small open neighborhood of
$([s_j])\in\prod_{j=1}^r\PP(I_{X,a_j})$ in {\it classical} topology.
\sk
We briefly explain the construction of $M_y(X)$. There is a flat family
$\Zc_T$ over a Zariski-open subset $T$ of $\C^r$ defined by
$$\Zc_T:=\Zc\cap(Y^{(0)}\times T)\q\h{with}\q
\Zc:=\mcap_{i=1}^r\{s_i=t_is'_i\}\subset Y^{(0)}\times\C^r.$$
Here $t_1,\dots,t_r$ are the coordinates of $\C^r$, and
$T\subset\C^r$ is defined by the following condition:
$t\in T\iff\dim\Zc_t=\dim Y^{(0)}-r$ with
$\Zc_t:=\Zc\cap(Y^{(0)}\times\{t\})$ for $t\in\C^r$.
(Note that $\Zc_T$ is a complete intersection in $Y^{(0)}\times T$, and
is flat over $T$, i.e. the $t_i-c_i$ for $i\in[1,r]$ form a regular
sequence in $\Oc_{\Zc_T,(y,c)}$ for any $(y,c)\in\Zc_T$. This follows
from a well-known theory of commutative algebra about regular sequences
and flat morphisms, see e.g. [Ei].)
\sk
We have $\Zc_0=Y^{(r)}=X$ by definition.
Applying the iterated nearby cycle functors, we define a mixed
Hodge module on $X$ by
$$\M'(s'_1,\dots,s'_r):=\psi_{t_r}\cdots\psi_{t_1}\Q_{h,\Zc_T}[\dim X],$$
which satisfies
$$T^{\,\vir}_{y*}(X)=(-1)^{\dim X}T_{y*}\bl(\M'(s'_1,\dots,s'_r)\br).$$
Indeed, the nearby cycle functor of mixed Hodge modules corresponds
to the Gysin morphism of Borel-Moore homology via the
transformation $\DR_y$ in (1.2.1), see Proposition~(3.3) below.
(Note that the latter transformation is denoted by $mC_y$ in [Sch2].)
We can moreover prove the injectivity of the canonical morphism of
mixed Hodge modules
$$\Q_{h,X}[\dim X]\into\M'(s'_1,\dots,s'_r),$$
together with the equality (0.2) at the level of the Grothendieck
group of mixed Hodge modules by increasing induction on $r$.
Then $\M(s'_1,\dots,s'_r)$ is defined to be the cokernel of the above
injection.
It may depend on the $s'_i\,(i\in[1,r])$ although its image by $T_{y*}$
(and hence $M_y(X)$) does not by using a one-parameter family together
with Lemma~(2.6) below, see also [Sch3], Cor.~3.7.
\ms
By [BSY] and (1.3.4) below, $T_{y*}(X)$ and $T^{\,\vir}_{y*}(X)$
respectively specialize for $y=-1$ to the MacPherson-Chern class
$c(X)$ (see [Ma]) and the virtual Chern class $c^{\vir}(X)$
(called the Fulton or Fulton-Johnson class, see [Fu2], [FJ]) with
rational coefficients. Specializing Theorem~1 and using [Sch1],
Prop.~5.21, we then get the following.
\ms\nin
{\bf Corollary~1.} {\it 
With the notation and the assumptions of Theorem~$1$, there is
$M(X)\in\HH_{\ssb}(\Si)$, called the Milnor class, and satisfying}
$$\aligned&c^{\vir}(X)-c(X)=(i_{\Si,X})_*M(X),\\
&M(X)=c\bl((i_{\Si\setminus X',\Si})_{!\,}\varphi_f
\Q_Y\br)+(i_{\Si',\Si})_*M(X')\raise4pt\h{}.\endaligned$$

Here $c\bl((i_{\Si\setminus X',\Si})_{!\,}\varphi_f\Q_Y\br)$
is the MacPherson-Chern class [Ma] with $\Q$-coefficients of the
associated constructible function. To get Corollary~1 with
$\Z$-coefficients (i.e.\ in $\HH_{\ssb}(\Si,\Z)$), we have to prove some
assertions in this paper (e.g.\ Proposition~(4.1)) with $\Q$ replaced
by $\Z$, and use [Ve2] (and [Ke1], [Ke2] if
$\HH_{\ssb}(X)=\CH_{\ssb}(X)$).
\ms
As special cases of Theorem~1, we also get the following.
\ms\nin
{\bf Corollary~2.} {\it 
With the notation and the assumptions of Theorem~$1$, assume
furthermore $r=1$, i.e. $Y$ is smooth, or $\dim\Si=0$, i.e. $X$
has only isolated singularities.
Then the second term in the right-hand side of $(0.2)$ vanishes so that
$$M_y(X)=T_{y*}\bl((i_{\Si\setminus X',\Si})_{!\,}\varphi_f
\Q_{h,Y}\br).\leqno(0.3)$$

If $\dim\Si=0$, then we can omit
$(i_{\Si\setminus X',\Si})_!$ $($since $\Si\cap X'=\emptyset)$, and
$$M_y(X)=\mopl_{x\in\Sing\,X}\,\chi_y(\widetilde{H}^{\ssb}(F_x)).
\leqno(0.4)$$
Here $\widetilde{H}^{\ssb}(F_x)$ denotes the reduced Milnor cohomology of
$f:(Y,x)\to(\C,0)$ endowed with a canonical mixed Hodge structure, and
$\chi_y$ is a polynomial defined by using its Hodge numbers, see $(1.2.7)$
below.}
\ms
Note that $\chi_y(\widetilde{H}^{\ssb}(F_x))$ is essentially the
Steenbrink spectrum [St3] with information on the action of the
monodromy forgotten (see [CMSS], Remark~3.7 for a more precise
statement).
This spectrum coincides with the spectrum of the complete intersection
$(Y^{(0)},X)$ at a general point of $C_XY^{(0)}$ over $x\in\Sing\,X$
defined in [DMS], if $s_1,\dots,s_{r-1}$ and $s'_1,\dots,s'_r$ are
sufficiently general.
One can calculate $\chi_y(\widetilde{H}^{\ssb}(F_x))$ by taking an
embedded resolution of singularities and using the theory of motivic
Milnor fibers [DL] (which is explained in a slightly different
situation in Proposition~(5.7) below). We can also apply a construction
of Steenbrink [St2] in the isolated singularity case or the theory of
mixed Hodge modules as in [BS], Th.~4.3.
\sk
By using [Sch1], Prop.~5.21, we see that Corollary~2 specializes for
$y=-1$ to a formula for the Chern classes in the hypersurface case
which was conjectured by S.~Yokura [Yo2] (see also [Yo3], [Yo4]), and
proved by A.~Parusi\'nski and P.~Pragacz [PP].
Note that the formulation itself of Theorem~0.2 in loc.~cit.\
cannot be generalized to the case of Hirzebruch classes because of
some monodromy problem.
The equality (0.4) specializes for $y=-1$ to a formula of J.~Seade and
T.~Suwa for the Chern classes in the case of complete intersections with
isolated singularities (see [SeSu], [Su]).
\ms
Note finally that the Hirzebruch-Milnor class $M_y(X)$ can be
{\it calculated explicitly} by using (0.2), see Section~5 below. This
is quite important form the view point of applications of Theorem~1.
\ms
Part of this work was done during a visit of the second named author
to Mathematical Institute of M\"unster University, and he thanks
the institute for the financial support.
The first named author is partially supported by NSF-1005338.
The second named author is partially supported by Kakenhi 21540037.
The third named author is supported by the SFB 878
``groups, geometry and actions''.
\ms
In Section 1 we review some basics of the Hirzebruch characteristic
classes. In Section 2 we prove some properties of the Gysin morphisms
which are used in later sections. In Section 3 we show some assertions
in the global complete intersection case, e.g. Proposition~(3.4).
In Section 4 we prove Theorem~1 after showing Proposition~(4.1).
In Section 5 we show how the Hirzebruch-Milnor classes can be
calculated explicitly.

\ms\nin
{\bf Conventions.}
\sk\nin
1. A variety means a complex algebraic variety which is always
assumed reduced and irreducible except for the case of
subvarieties (e.g. hypersurfaces) which may be non-reduced or reducible.
We assume a variety is quasi-projective when we use [Fu1].
\sk\nin
2. Cohomology and homology groups are always with $\Q$-coefficients
unless the coefficients are explicitly stated.
\sk\nin
3. In this paper, $\HH_k(X)$, $\HH^k(X)$ are respectively
$H^{\BM}_{2k}(X,\Q)$ and $H^{2k}(X,\Q)$ in order to simplify the
explanations. We have similar assertions for $\HH_k(X):=\CH_k(X)_{\Q}$
where $\HH^k(X)$ is either the Chow cohomology in [Fu1] or
Fulton-MacPherson's operational Chow cohomology in [FM], [Fu2]
(see [To] for the relation with $H^{2k}(X,\Q)$).
\sk\nin
4. The nearby and vanishing cycle functors $\psi_f$, $\varphi_f$
for $D^b\MHM(X)$ are not shifted by $-1$ as in [Sa1], [Sa2].
They are compatible with the corresponding functors for the
underlying $\Q$-complexes without a shift of complexes, but do not
preserve mixed Hodge modules.
\sk\nin
5. We use {\it left} $\Dc$-modules, although right $\Dc$-modules are
mainly used in [Sa1], [Sa2].
The transformation between filtered left and right $\Dc$-modules is
given by
$$(M,F)\mapsto(\Omega^{\dim X}_X,F)\otimes_{\Oc_X}(M,F),$$
where $(M,F)$ are left $\Dc_X$-modules, and $\Gr_p^F\Omega_X^{\dim X}=0$
for $p\ne-\dim X$.

\bs\bs
\centerline{\bf 1. Hirzebruch characteristic classes}
\bs\nin
In this section we review some basics of the Hirzebruch characteristic
classes.
\ms\nin
{\bf 1.1.~Cohomology Hirzebruch classes.}
Let $X$ be a smooth complex algebraic variety of dimension $n$.
The cohomology Hirzebruch characteristic class $T_y^*(TX)$ of
the tangent bundle $TX$ is defined by
$$T^*_y(TX):=\mprod_{i=1}^n\,Q_y(\al_i)\in \HH^{\ssb}(X)[y],
\leqno(1.1.1)$$
where $Q_y(\al)$ is explained as below and
the $\{\al_i\}$ are the (formal) Chern roots of $TX$, i.e.
$$\mprod_{i=1}^n(1+\al_it)=\msum_{j=0}^n\,c_j(TX)t^j.$$

We have normalized and unnormalized power series
(see [Hi], 1.8 and [HBJ], 5.4.):
$$Q_y(\al):=\frac{\al(1+y)}{1-e^{-\al(1+y)}}-\al y,\,\,\,
\Qt_y(\al):=\frac{\al(1+ye^{-\al})}{1-e^{-\al}}\in
\Q[y][[\al]].
\leqno(1.1.2)$$
Their initial terms are as follows:
$$Q_y(0)=1,\q\Qt_y(0)=1+y.
\leqno(1.1.3)$$
These two power series have the following relation
$$Q_y(\al)=(1+y)^{-1}\,\Qt_y(\al(1+y)).
\leqno(1.1.4)$$

\ms\nin
{\bf 1.2.~Homology Hirzebruch classes.}
Let $X$ be a complex algebraic variety.
Let $\MHM(X)$ be the category of mixed Hodge modules on $X$.
For $\MM\in D^b\MHM(X)$, its homology Hirzebruch characteristic
class is defined by
$$\aligned T_{y*}(\MM)&:=td_{(1+y)*}\bl(\DR_y[\MM]\br)\in
\HH_{\ssb}(X)\bl[y,\h{$\frac{1}{y(y+1)}$}\br]\q\h{with}\\
\DR_y[\MM]&:=\msum_{i,p}\,(-1)^i\,\bl[\Hc^i\Gr_F^p\DR(\MM)\br]
\,(-y)^p\in K_0(X)[y,y^{-1}],\raise12pt\h{ }\endaligned
\leqno(1.2.1)$$
setting $F^p:=F_{-p}$. Here we define
$$td_{(1+y)*}:K_0(X)[y,y^{-1}]\to
\HH_{\ssb}(X)\bl[y,\h{$\frac{1}{y(y+1)}$}\br]$$
to be the scalar extension of the Todd class transformation
$$td_*:K_0(X)\to\HH_{\ssb}(X)$$
(denoted by $\tau$ in [BFM]) together with the multiplication by
$(1+y)^{-k}$ on the degree $k$ part (see [BSY]).
The last multiplication is closely related to the identity (1.1.4), and
we have actually by [Sch1], Prop.~5.21
$$T_{y*}(\MM)\in \HH_{\ssb}(X)[y,y^{-1}].$$
In this paper, we set $\HH_k(X):=H^{\BM}_{2k}(X,\Q)$ to simplify the
arguments, see Convention~3.
In some other papers (as [BSY], [CMSS], etc.), $\DR_y[\MM]$ is denoted
by $mC_y(\MM)$.
(Note also that it is nontrivial to show that the $\Hc^i\Gr_F^p\DR(\MM)$
are $\Oc_X$-modules in the singular case, see [Sa1], Lemma~3.2.6.)

By definition we have for $k\in\Z$
$$\DR_y[\MM(k)]=\DR_y[\MM]\,(-y)^{-k},\q
T_{y*}(\MM(k))=T_{y*}(\MM)\,(-y)^{-k},
\leqno(1.2.2)$$
where the Tate twist $(k)$ on the filtered $\Dc$-module part is given
by the shift of the filtration $[k]$ which is defined by
$$(F[k])^p:=F^{p+k},\q(F[k])_p:=F_{p-k}.
\leqno(1.2.3)$$

The {\it homology Hirzebruch characteristic class} $T_{y*}(X)$ is
defined by applying the above definition to the case $\MM=\Q_{h,X}$
(see [BSY]), i.e.
$$\aligned &T_{y*}(X):=T_{y*}(\Q_{h,X})=td_{(1+y)*}\DR_y[X]\in
\HH_{\ssb}(X)[y],\\
&\h{with}\q\DR_y[X]:=\DR_y[\Q_{h,X}].\endaligned$$
This coincides with the definition using the Du Bois complex [DB] by
[Sa3], and it is known that $T_{y*}(X)$ belongs to $\HH_{\ssb}(X)[y]$,
see [BSY]. In case $X$ is smooth, we have
$$\DR_y[X]=\La_y[T^*X],
\leqno(1.2.4)$$
where $\La_y[V]:=\sum_p[\La^pV]\,y^p$ for a vector bundle $V$.
Indeed, we have
$$\DR(\Q_{h,X})=\DR(\Oc_X)[-n]=\Omega_X^{\ssb}\q\h{with}\,\,\,
n:=\dim X,
\leqno(1.2.5)$$
and the Hodge filtration $F^p$ on $\Omega_X^{\ssb}$ is given by the
truncation $\si_{\ge p}$ in [D1].

To show the coincidence with the above definition of $T_{y*}(X)$ in
the smooth case, we need the relation (1.1.4) together with some
calculation about Hirzebruch's power series $Q_y(\al)$ as in [HBJ],
Sect.~5.4 or in the proof of [Yo1], Lemma~2.3.7, which is closely
related to the generalized Hirzebruch-Riemann-Roch theorem in [Hi],
Th.~21.3.1.
This coincidence is part of the characterization of the Hirzebruch
characteristic classes for singular varieties, see [BSY], Th.~3.1.

By a generalization of the Riemann-Roch theorem [BFM], 
$td_*$ commutes with the pushforward under proper morphisms, and
so does $T_{y*}$.
If we apply this to $a_X:X\to pt$ with $X$ compact, then the
pushforward for $\HH_{\ssb}$ is identified with the degree map or the
trace morphism, and we have
$$K_0(pt)=\Z,\q\HH_{\ssb}(pt)=\Q,\q\MHM(pt)=\MHS,
\leqno(1.2.6)$$
where $\MHS$ is the category of graded-polarizable mixed $\Q$-Hodge
structures in [D1].
By definition we have for $H^{\ssb}\in D^b\MHS$
$$T_{y*}(H^{\ssb})=\chi_y(H^{\ssb}):=
\msum_{j,p}\,(-1)^j\dim_{\C}\Gr_F^pH_{\C}^j\,(-y)^p.
\leqno(1.2.7)$$
So the degree 0 part of $T_{y*}(X)$ is identified with
$\chi_y\bl(H^{\ssb}(X)\br)$ if $X$ is compact and connected.
Here the cohomology groups $H^{\ssb}(X)$ (with the zero differential)
is defined in $D^b\MHS$ by [D1], and this is compatible with [Sa2] by
[Sa3] (using [Ca]).

\ms\nin
{\bf 1.3.~Virtual Hirzebruch classes.}
Hirzebruch also introduced the notion of virtual $\chi_y$-genus
(or $\chi_y$-characteristic) which is the $\chi_y$-genus of general
complete intersections $X$ in smooth projective varieties $Y$.
Let $X$ be a complete intersection in a smooth projective variety $Y$.
The virtual Hirzebruch characteristic class $T^{\,\vir}_{y*}(X)$ can
be defined like the virtual genus by
$$T^{\,\vir}_{y*}(X):=td_{(1+y)*}\DR^{\vir}_y[X]\in\HH_{\ssb}(X)[y],
\leqno(1.3.1)$$
with $\DR^{\vir}_y[X]$ the image in $K_0(X)[[y]]$ of 
$$\La_y(T^*_{\vir}X)=\La_y[T^*Y|_X]/\La_y[N^*_{X/Y}]\in K^0(X)[[y]],
\leqno(1.3.2)$$
and $\DR^{\vir}_y[X]$ belongs to $K_0(X)[y]$ by Proposition~(3.4)
below.
Here $K^0(X)$, $K_0(X)$ are respectively the Grothendieck group of
locally free sheaves of finite length and that of coherent sheaves.
We denote respectively by $T^*Y$ and $N^*_{X/Y}$ the cotangent and
conormal bundles, and the virtual cotangent bundle is defined by
$$T^*_{\vir}X:=[T^*Y|_X]-[N^*_{X/Y}]\in K^0(X).$$
More precisely, $N^*_{X/Y}$ in the non-reduced case is defined by the
locally free sheaf $\I_X/\I_X^2$ on $X$ where $\I_X\subset\Oc_Y$ is
the ideal sheaf of the subvariety $X$ of $Y$.
Here we set for a virtual vector bundle $V$ on $X$ in general
$$\La_yV:=\msum_{p\ge 0}\,\La^pV\,y^p\in K^0(X)[[y]].
\leqno(1.3.3)$$

We can also define $T^{\,\vir}_{y*}(X)$ by using the virtual
tangent bundle
$$T_{\vir}X:=[TY|_X]-[N_{X/Y}]\in K^0(X),$$
together with the above cohomological transformation $T^*_y$ as in
[CMSS] (in the hypersurface case) so that
$$T^{\,\vir}_{y*}(X)=T^*_y(T_{\vir}X)\cap[X]\q\h{in}\,\,\,
\HH_{\ssb}(X)[y],
\leqno(1.3.4)$$
see Proposition~(1.4) below.
Here $[X]:=\sum_jm_j[X_j]$ is the fundamental class of $X$ with $m_j$
the multiplicities along the irreducible components $X_j$ of $X$.
We have the equality $T_{y*}(X)=T^{\,\vir}_{y*}(X)$ if $X$ is
smooth.
The problem is then how to describe the difference in the singular
case, which is called the {\it Hirzebruch-Milnor class}.
(For the degree-zero part, i.e.\ on the level of Hodge polynomials,
see also [LM].)
\ms\nin
{\bf 1.4.~Proposition.} {\it The equality $(1.3.4)$ holds in
$\HH_{\ssb}(X)[[y]]$ under the assumption that $X$ is a locally complete
intersection in a smooth variety $Y$ of pure codimension $r$.}
\ms\nin
{\it Proof.}
By the compatibility of $td_*$ with the cap product [BFM], we have
$$td_*\bl(\DR^{\vir}_y[X]\br)=ch(\Lambda_yT^*_{\vir}X)\cap td_*(X)
\q\h{in}\,\,\,\HH_{\ssb}(X)[[y]],$$
where $\DR^{\vir}_y[X]$ is as in (1.3.1), and $ch$ denotes also the
scalar extension of the Chern character $ch:K^0(X)\to\HH^{\ssb}(X)$
under $\Q\into\Q[[y]]$. By [Ve1], Th.~7.1, we have
$$td_*(X):=td_*([\Oc_X])=td_*(i^*[\Oc_Y])=
i^!td_*(Y)\cap td^*(N_{X/Y})^{-1},$$
and
$$td_*(Y)=td^*(TY)\cap[Y].$$
Using the relation of $i^!$, $i^*$ with the cap product (see (2.3.4)
below), we then get
$$td_*\bl(\DR^{\vir}_y[X]\br)=
ch(\Lambda_yT^*_{\vir}X)\cdot td^*(T_{\vir}X)\cap[X].$$

Let $\al_i$ and $\be_j$ be respectively the (formal) Chern roots of
$TY|_X$ and $N_{X/Y}$. Then
$$ch(\Lambda_yT^*_{\vir}X)\cdot td^*(T_{\vir}X)=
\frac{\mprod_{i=1}^m(1+ye^{-\al_i})\cdot\al_i/(1-e^{-\al_i})}
{\mprod_{j=1}^r(1+ye^{-\be_i})\cdot\be_j/(1-e^{-\be_j})},$$
where $m=\dim Y$ and $r=m-n$ with $n=\dim X$.
Indeed, if $\ell(\gamma_j)$ denotes a (formal) line bundle class with the
first Chern class $\gamma_j$ in some ring extension of $K^0(X)$, then
$$ch\bl(\ell(\gamma_1)\,\cdots\,\ell(\gamma_k)\br)=
ch\bl(\ell(\gamma_1+\cdots+\gamma_k)\br)=e^{\gamma_1+\cdots+\gamma_k},$$
and this implies (see also [Hi])
$$ch(\Lambda_yT^*_{\vir}X)=
\mprod_{i=1}^m(1+ye^{-\al_i})/\mprod_{j=1}^r(1+ye^{-\be_i}).$$

Since $td_{(1+y)*}$ is the composition of $td_*$ with the multiplication
by $(1+y)^{-k}$ on $\HH_k(X)$, and the last multiplication corresponds
to the multiplication by $(1+y)^k$ on $\HH^k(X)$ under the cap product,
we then get the desired equality as follows:
$$\aligned &\q\,\,\,td_{(1+y)*}\bl(\DR^{\vir}_y[X]\br)\\
&=\frac{\mprod_{i=1}^m\bl(1+ye^{-\al_i(1+y)}\br)\cdot
\al_i(1+y)/(1-e^{-\al_i(1+y)})}
{\mprod_{j=1}^r\bl(1+ye^{-\be_i(1+y)}\br)\cdot
\be_j(1+y)/(1-e^{-\be_j(1+y)})}\cap(1+y)^{-n}[X]\\
&=\frac{\mprod_{i=1}^m\bl(\Qt_y(\al_i(1+y))/(1+y)\br)}
{\mprod_{j=1}^r\bl(\Qt_y(\be_j(1+y))/(1+y)\br)}\cap[X]\\
&=T^*_y\bl([TY|_X]-[N_{X/Y}]\br)\cap[X]\\
&=T^*_y(T_{\vir}X)\cap[X],\endaligned$$
where the relation (1.1.4) between the power series $Q_y$ and $\Qt_y$ is
used, see also [BSY], (1.1) and [Yo1], Lemma~2.3.7.
This finishes the proof of Proposition~(1.4).

\bs\bs
\centerline{\bf 2. Gysin morphisms}
\bs\nin
In this section we prove some properties of the Gysin morphisms which
are used in later sections.
\ms\nin
{\bf 2.1.~Construction.}
Let $i:X\into Y$ be a locally complete intersection morphism of pure
codimension $r$. By [Fu2], [Ve1] there are Gysin morphisms
$$i^!:\CH_k(Y)\to\CH_{k-r}(X),\q i^!:H^{\BM}_{2k}(Y)\to
H^{\BM}_{2k-2r}(X)(r),
\leqno(2.1.1)$$
in a compatible way with the cycle map.
These are constructed by using the deformation to the normal cone
$q:Z\to S:=\C$ such that $q^{-1}(0)=N_{X/Y}$ and
$q^{-1}(S')=Y\mtim S'$ where $S':=\C^*$.
More precisely, $Z$ is the complement of the proper transform of
$Y\mtim\{0\}$ in the blow-up of $Y\mtim S$ along $X\mtim\{0\}$,
see loc.~cit.

We recall here the definition of the Gysin morphism for Borel-Moore
homology (see [Ve1]).
We have the isomorphisms
$$(R^kq_{!\,}\Q_Z)_0=H^k_c(N_{X/Y}),\q
R^kq_{!\,}\Q_Z|_{S'}=a_{S'}^{-1}H^k_c(Y),$$
where $a_{S'}:S'\to pt$ is the natural morphism.
So we get the cospecialization morphism
$$\spe^*:H^k_c(N_{X/Y})\to H^k_c(Y).
\leqno(2.1.3)$$
Taking the dual, we get the Gysin morphism
$$i^!:H^{\BM}_k(Y)\buildrel{\spe_*}\over\longrightarrow
H^{\BM}_k(N_{X/Y})\buildrel{\sim}\over\leftarrow H^{\BM}_{k-2r}(X)(r),
\leqno(2.1.4)$$
where the last isomorphism is induced by the pullback under the
canonical morphism $\pi_N:N_{X/Y}\to X$.

The above definition implies the compatibility with the
pushforward under a proper morphism $\rho:\Yt\to Y$ such that
$\X=X\mtim_Y\Yt$ has codimension $r$ everywhere in $\Yt$.
In this case we have the cartesian diagram
$$\begin{matrix}\X&\buildrel{\iti}\over\into&\Yt\\
\,\,\,\,\downarrow\!{\scriptstyle\rho'}&&
\,\,\,\downarrow\!{\scriptstyle\rho}\\
X&\buildrel{i}\over\into&Y\end{matrix}$$
and
$$i^!\scc\rho_*=\rho'_*\scc\iti^{\,!}.
\leqno(2.1.5)$$
Indeed, the specialization morphism $\spe_*$ in (2.1.4) is compatible
with the pushforward by proper morphisms, and $N_{\X/\Yt}\to\X$ is the
base change of $N_{X/Y}\to X$ by the hypothesis.
\ms\nin
{\bf 2.2.~Sheaf-theoretic description} [Ve1].
With the above notation, set $Z_s:=q^{-1}(s)$, and let
$i_s:Z_s\into Z$ denote the inclusion for $s\in S$.
We have
$$Z_0=N_{X/Y},\q Z_s=Y\,(s\ne 0).
\leqno(2.2.1)$$
Let $p:Z\to Y$ be the canonical projection.
Set $p_s:=p\scc i_s:Z_s\to Y$.

For $s\in S$ we have the canonical flasque resolutions
$$\Q_Z\buildrel{\rm qi}\over\longrightarrow\I^{\ssb}_Z,\q
\Q_{Z_s}\buildrel{\rm qi}\over\longrightarrow\I^{\ssb}_{Z_s},$$
together with the canonical morphisms
$$\I^{\ssb}_Z\to(i_s)_*\I^{\ssb}_{Z_s},
\leqno(2.2.2)$$
which are compatible with $\Q_Z\to(i_s)_*\Q_{Z_s}$.
Indeed, $(i_s)_*\I^{\ssb}_{Z_s}$ is naturally isomorphic to the
canonical flasque resolution of $(i_s)_*\Q_{Z_s}$, and this resolution
is functorial for morphisms of sheaves.
The dualizing complex $\D^{\ssb}_Z$ is defined by
$$\D^k_Z(V):=\Hom(\Ga_c(V,\I^{-k}_Z),\Q)\q\h{for open}\,\,\,
V\subset S,$$
and similarly for $\D^{\ssb}_{Z_s}$.
Then we have by (2.2.2) the canonical morphisms
$$(i_s)_*\D^{\ssb}_{Z_s}\to\D^{\ssb}_Z.$$
They induce the canonical morphisms
$$\Ga(Z_s,\D^{\ssb}_{Z_s})\to\Ga_c(S,q_*\D^{\ssb}_Z),
\leqno(2.2.3)$$
and this is a quasi-isomorphism for $s=0$.
We thus get the canonical morphisms of complexes
$$\Ga(Y,\D^{\ssb}_Y)=\Ga(Z_1,\D^{\ssb}_{Z_1})\to
\Ga_c(S,q_*\D^{\ssb}_Z)]\buildrel{\rm qi}\over\longleftarrow
\Ga(Z_0,\D^{\ssb}_{Z_0})=\Ga(N_{X/Y},\D^{\ssb}_{N_{X/Y}}),
\leqno(2.2.4)$$
where qi means an quasi-isomorphism.

For any open subset $U\subset Y$, set
$$X_U:=X\cap U,\q Z_U:=p^{-1}(U),\q Z_{U,s}:=Z_U\cap Z_s,\q
q_U:=q|_{Z_U}:Z_U\to S.$$
Then
$$q_U^{-1}(0)=N_{X_U/U},\q q_U^{-1}(s)=U\,\,\,(s\ne 0),$$
and (2.2.4) implies the canonical morphisms of complexes
$$\Ga(U,\D^{\ssb}_U)\to\Ga_c(S,(q_U)_*\D^{\ssb}_{Z_U})
\buildrel{\rm qi}\over\longleftarrow
\Ga(N_{X_U/U},\D^{\ssb}_{N_{X_U/U}}).
\leqno(2.2.5)$$
Here we have presheaves defined by associating to $U$ respectively
$$\Ga_c(S,(q_U)_*\D^k_{Z_U}),\q\Ga(N_{X_U/U},\D^k_{N_{X_U/U}}),$$
and these are sheaves.
So we get the canonical morphisms of sheaf complexes on $Y$,
and this gives the morphism in the derived category $D^b_c(Y,\Q_Y)$
$$\D^{\ssb}_Y\to(i\scc\pi_N)_*\D^{\ssb}_{N_{X/Y}}
\buildrel{\sim}\over\leftarrow i_*\D^{\ssb}_X(r)[2r],
\leqno(2.2.6)$$
where $\pi_N:N_{X/Y}\to X$ denotes the canonical projection
and the last quasi-isomorphism is given by $\pi_N^*$.
So we get the sheaf-theoretic description of the Gysin morphism in (2.1).
(This is a detailed version of a slight modification of a sketch of an
argument in [Ve1], Sect.~8 where a point of $S'=\C\setminus\{0\}$ was
not chosen as in the above argument.
We may assume that this point is $1$ by using an equivariant
$\C^*$-action on $q:Z\to S$.)
\ms\nin
{\bf 2.3.~Compatibility with the cap product.}
Let $\zeta\in H^j(Y,\Q)$.
By the canonical isomorphism
$$H^j(Y,\Q)=\Hom(\Q_Y,\Q_Y[j]),
\leqno(2.3.1)$$
together with $\D_Y=\Q_Y\otimes\D_Y$, we get the morphism
$$\zeta\,\cap:H^{\BM}_k(Y,\Q)\to H^{\BM}_{k-j}(Y,\Q).
\leqno(2.3.2)$$
So the above sheaf-theoretic construction of the Gysin morphism implies
the compatibility with the cap product
$$\begin{matrix}H^{\BM}_k(Y,\Q)&\buildrel{i^!}\over\to&
H^{\BM}_{k-2r}(X,\Q)(r)\\
\q\downarrow{\!\scriptstyle{\zeta\,\cap}}&&
\q\downarrow{\!\scriptstyle{i^*\zeta\,\cap}}\\
H^{\BM}_{k-j}(Y)&\buildrel{i^!}\over\to&H^{\BM}_{k-j-2r}(X)(r)\end{matrix}
\leqno(2.3.3)$$
i.e.
$$i^!(\zeta\cap\xi)=i^*\zeta\cap i^!\xi\q\h{for}\,\,\,\zeta\in
H^j(Y,\Q),\,\,\xi\in H^{\BM}_k(Y).
\leqno(2.3.4)$$
Since $i^![Y]=[X]$, this implies the commutative diagram
$$\begin{matrix}H^{2n-k}(Y)(n)&\buildrel{i^*}\over\to&H^{2n-k}(X)(n)\\
\q\downarrow{\!\scriptstyle{\cap[Y]}}&&
\q\downarrow{\!\scriptstyle{\cap[X]}}\\
H^{\BM}_k(Y)&\buildrel{i^!}\over\to&H^{\BM}_{k-2r}(X)(r)\end{matrix}
\leqno(2.3.5)$$
This is used in the proof of [CMSS], Lemma~3.1 in the hypersurface
case (i.e. $r=1$).

\ms\nin
{\bf 2.4.~Principal divisor case.}
Assume $X$ is a globally principal divisor on $Y$, i.e. there is a
function on $Y$ with $X=f^{-1}(0)$ and $r=1$.
Then the Gysin morphism induces a `well-defined' morphism
$$i^!:\CH_k(Y\setminus X)\to\CH_{k-1}(X),
\leqno(2.4.1)$$
using the exact sequence
$$\CH_k(X)\buildrel{i_*}\over\to\CH_k(Y)\buildrel{j^*}\over\to
\CH_k(Y\setminus X)\to 0,$$
where `well-defined' means that $i^!:\CH_k(Y)\to\CH_{k-1}(X)$ factors
through the surjection $j^*$.
Indeed, the composition $i^!\scc i_*$ vanishes by the triviality
of the normal bundle.

By [Ve1], Th.~7.1, we have
$$td_*\scc i^*=i^!\scc td_*.
\leqno(2.4.2)$$
where $i^*:K_0(Y)\to K_0(X)$ is the pull-back of the Grothendieck
groups of coherent sheaves defined by the mapping cone of the
multiplication by the global defining function $f$ of $X$.
Using the exact sequence
$$K_0(X)\buildrel{i_*}\over\to K_0(Y)\buildrel{j^*}\over\to
K_0(Y\setminus X)\to 0,$$
the last definition implies that the composition $i^*\scc i_*$ vanishes
so that $i^*$ induces also a `well-defined' morphism
$$i^*:K_0(Y\setminus X)\to K_0(X).
\leqno(2.4.3)$$

We have a similar property for the nearby cycle
functor $\psi_f$ with $f$ a global defining function of $X$, and we
have the `well-defined' functor
$$\psi_f[-1]:\MHM(Y\setminus X)\to\MHM(X).
\leqno(2.4.4)$$
Indeed, the restriction gives the surjection
$\MHM(Y)\to\MHM(Y\setminus X)$ by extendability of mixed Hodge modules
in the algebraic case, and $\psi_f\M'[-1]$ is independent of the
extension $\M'$ of $\M$ to $Y$ by using the functorial morphism
$id\to j_*j^*$, since $\psi_fM''[-1]$ vanishes if
${\rm supp}\,M''\subset X$.

As for Borel-Moore cohomology, we have the long exact sequence
$$\cdots\to H^{\BM}_k(X)\buildrel{i_*}\over\to H^{\BM}_k(Y)
\buildrel{j^*}\over\to H^{\BM}_k(Y\setminus X)\to H^{\BM}_{k-1}(X)
\to\cdots,
\leqno(2.4.5)$$
and the next proposition shows the vanishing of the composition
$$H^{\BM}_k(X)\buildrel{i_*}\over\to H^{\BM}_k(Y)
\buildrel{i^!}\over\to H^{\BM}_{k-2}(X)(1).
\leqno(2.4.6)$$
\ms\nin
{\bf 2.5.~Proposition.} {\it With the above notation and the assumptions,
the composition of the two morphisms in $(2.4.6)$ vanishes.}
\ms\nin
{\it Proof.} By the construction in (2.1), the restriction of $\spe_*$
to the image of $i_*$ is essentially the identity.
More precisely, we have
$$\spe_*\scc i_*=(s_0)_*:H^{\BM}_k(X)\to H^{\BM}_k(N_{X/Y}),
\leqno(2.5.1)$$
where $s_0:X\into N_{X/Y}$ is the zero-section.
We have moreover
$$s_0^!=(\pi_N^*)^{-1}:H^{\BM}_k(N_{X/Y})\simto H^{\BM}_{k-2}(X)(1),
\leqno(2.5.2)$$
where $\pi_N:N_{X/Y}\to X$ is the projection.
Indeed, we have
$$s_0^!\scc \pi_N^*=id:H^{\BM}_k(X)\to H^{\BM}_k(X),
\leqno(2.5.3)$$
since the deformation to the normal cone is identified with a
trivial deformation.

So the assertion is reduced to
$$s_0^!\scc(s_0)_*=0.
\leqno(2.5.4)$$
Since the normal bundle $N_{X/Y}$ is trivial, we have a section
$s_1$ whose image is disjoint from that of $s_0$, and
$s_1^!\scc \pi_N^*=id$ by the same argument as above.
So the assertion is further reduced to
$$s_1^!\scc(s_0)_*=0.
\leqno(2.5.5)$$
This is more or less trivial by the construction of $s_1^!$
since the images of $s_0$ and $s_1$ are disjoint.
(The reader can also use (2.1.5) together with the fact that the
codimension of an empty set can be any number.)
This finishes the proof of Proposition~(2.5).
\ms
By a similar argument, we have the following.
\ms\nin
{\bf 2.6.~Lemma.} {\it Let $X$ be a complex algebraic variety.
Let $i_t:X=X\times\{t\}\into X\times\C$ denote the inclusion for
$t\in\C$.
Then $i_t^!:\HH_{\ssb}(X\times\C)\to\HH_{\ssb-1}(X)$ is independent of
$t\in\C$.}
\ms\nin
{\it Proof.} Let $Y:=X\times\C$ so that $Y=N_{X/Y}$ where
$X=X\times\{t\}\subset Y$. Then the assertion follows from (2.5.2). 
\bs\bs
\centerline{\bf 3. Global complete intersections}
\bs\nin
In this section we show some assertions in the global complete
intersection case, e.g. Proposition~(3.4).
We first give proofs of the following lemma and proposition for the
convenience of the reader. These are related to the independence of
the choice of the sections $s_i$.
They should be well-known to specialists.
\ms\nin
{\bf 3.1.~Lemma.} {\it Let $L$ be a line bundle on a smooth
variety $X$, and $s_i\in\Ga(X,L)$ for $i\in[1,m]$.
Assume $\bigcap_{i\in[1,m]}s_i^{-1}(0)=\emptyset$.
Then there is a non-empty Zariski-open subset $U$ of $\PP^{m-1}$
such that the zero locus of $\sum_{i\in[1,m]}t_is_i$ in $X$ is
smooth for any $t=(t_1,\dots,t_m)\in U$.}
\ms\nin
{\it Proof.} Consider
$$Z:=\bl\{(y,t)\in X\times\PP^{m-1}\,\big|\,
\msum_{i\in[1,m]}\,t_is_i(y)=0\br\}.$$
By the Sard-type theorem, it is enough to show that $Z$ is smooth.
Set
$$U_j:=\{t_j\ne 0\}\subset\PP^{m-1},\q
X_i:=X\setminus s_i^{-1}(0).$$
Then $X=\bigcup_iX_i$, and
$Z\cap\bl(X_i\mtim\PP^{m-1}\br)$ is covered by $X_i\times U_j$
with $i\ne j$ by the definition of $Z$.
Moreover, $Z\cap\bl(X_i\times U_j\br)$ is defined by
$$t_i/t_j=-\msum_{k\ne i}\,(t_k/t_j)(s_k(y)/s_i(y))
\q\h{in}\,\,\,X_i\times U_j,$$
where the $t_k/t_j\,(k\ne j)$ are the affine coordinates of $U_j$.
So the assertion follows.
\ms\nin
{\bf 3.2.~Proposition.} {\it
Let $L$ be a very ample line bundle on a projective variety $Y$.
Let $s_i\in\Ga(Y,L^{\otimes a_i})\,(i=1,\dots,r)$ with $a_i$
a decreasing sequence of positive integers where $r\ge 2$.
Set $Y^{(j)}:=\bigcap_{i\in[1,j]}s_i^{-1}(0)\,(j=1,\dots,r)$.
Assume $Y\setminus Y^{(r)}$ is smooth and ${\rm codim}_YY^{(r)}=r$.
Then, the $Y^{(j)}\setminus Y^{(r)}$ are smooth for any $j\in[1,r-1]$
by replacing $s_i$ with $s_i+\msum_{i'>i}\,s'_{i,i'}s_{i'}$
if we choose the $s'_{i,i'}\in\Ga(Y,L^{\otimes(a_i-a_{i'})})$
generically.}
\ms\nin
{\it Proof.} We proceed by increasing induction on $r\ge 2$.
Consider the embedding $Y\into\PP^N$ by the line bundle $L$. Let
$\{\tau^{(i)}_k\}_{k\in[1,\nu_i]}$ be a basis of $\Ga(Y,L^{\otimes i})$.
We apply Lemma~(3.1) to the case where $X=Y\setminus Y^{(r)}$,
$m=\sum_{i=1}^r\nu_{a_1-a_i}$, and the $s_j\,(j\in[1,m])$ in Lemma~(3.1)
are given by the restrictions of
$$\tau^{(a_1-a_i)}_ks_i\q(k\in[1,\nu_{a_1-a_i}],\,i\in[1,r]).$$
Note that the intersection of the zero loci of
$\tau^{(a_1-a_i)}_ks_i\,(k\in[1,\nu_{a_1-a_i}])$ coincides with
$s_i^{-1}(0)$ for each $i\in[1,r]$.
We may assume that the coefficient of $s_1$ does not vanish
since this condition defines a non-empty Zariski-open subset of
$\PP^{m-1}$ (and the intersection of any two non-empty Zariski-open
subsets of $\PP^{m-1}$ is non-empty).
So we get the assertion for $j=1$, and hence for $r=2$.
Moreover, the assertion for $r\ge 3,j\ge 2$ is reduced to the
assertion with $Y$ replaced by $Y^{(1)}$, and $(r,j)$ by $(r-1,j-1)$.
Thus we can proceed by induction on $r$.
This finishes the proof of Proposition~(3.2).

\ms
The following was shown in [Sch2] without assuming the condition that
the variety $Y$ is embeddable into a smooth variety.
We give here a simplified proof assuming this condition.
This proposition will be used in the proof of Proposition~(3.4) below.
\ms\nin
{\bf 3.3.~Proposition}~([Sch2]).
{\it Let $f:Y\to\C$ be a non-constant function on a complex algebraic
variety. Assume $Y$ is embeddable into a smooth variety.
Set $Y_0:=f^{-1}(0)$ with the inclusion $i_0:Y_0\into Y$. Then, for any
bounded complex of mixed Hodge modules $\M^{\ssb}$ on $Y$, we have
$$\aligned i_0^*\DR_y[\M^{\ssb}]&=\DR_y[\psi_f\M^{\ssb}](1+y)
\q\h{in}\,\,\,K_0(Y_0)[y,y^{-1}],\\
i_0^!T_{y*}(\M^{\ssb})&=T_{y*}(\psi_f\M^{\ssb})
\q\h{in}\,\,\,\HH_{\ssb}(Y_0)[y,y^{-1}],\endaligned
\leqno(3.3.1)$$
where $\psi_f\M^{\ssb}$ is viewed as a complex of mixed Hodge modules on
$Y_0$.}
\ms\nin
{\it Proof.} It is enough to show the first equality since the second
equality follows from it using [Ve1], Th.~7.1. (Here the term $(1+y)$
disappears since $i_0^!$ send $\HH_k(Y)$ to $\HH_{k-1}(Y_0)$.)
We may assume that $\M^{\ssb}$ is a mixed Hodge module $\M$ on $Y$. By
assumption there is a closed embedding $Y\into Y'$ with $Y'$ smooth. Set
$$Z:=Y\times\C\subset Z':=Y'\times\C.$$
Let $i_f:Y\into Z$ be the graph embedding, and $(M,F)$ be the underlying
filtered left $\Dc_{Z'}$-module of the direct image of $\M$ by $i_f$,
see [Sa1]. It has an increasing filtration $V$ of Kashiwara and
Malgrange such that $\partial_tt-a$ is nilpotent on $\Gr^a_VM$.
(In this paper we use left $\Dc$-modules, and $V^a$ corresponds to
$V_{-a}$ in loc.~cit., see Convention~5.)

Setting $Z'{}^*:=Y'\mtim\C^*$, we have $V^aM|_{Z'{}^*}=M|_{Z'{}^*}$ for
any $a\in\Q$. (Indeed, for any $m\in M$, we have $t^im\in V^aM$ for
$i\gg a$ and $t$ is invertible on $Z'{}^*$).
Using the well-definedness of (2.4.3), the left-hand side of (3.3.1)
can be given by applying $i_0^*$ to the mapping cone
$$C\bl(\Gr_F^{\ssb}\DR_{Z'/\C}(V^aM)\buildrel{\Gr\,\partial_t}\over\too
\Gr_F^{\ssb-1}\DR_{Z'/\C}(V^bM)\br),
\leqno(3.3.2)$$
and considering it in $K_0(Y_0)[y,y^{-1}]$ (using [Sa1], Lemma~3.2.6).
Here $a,b$ can be any rational numbers satisfying
$$a-1\ge b>0.
\leqno(3.3.3)$$
(This implies that [Sa1], (3.2.1.3) is not needed in [Sch2].)
By [Sa1], (3.2.1.2), we have the isomorphisms
$$t:V^aM\simto V^{a+1}M\q\h{for any}\,\,\,a>0.
\leqno(3.3.4)$$
So the left-hand side of (3.3.1) is given by
$$-\DR_y[(V^a/V^{a+1})(M,F)]+\DR_y[(V^b/V^{b+1})(M,F[-1])],
\leqno(3.3.5)$$
where $(V^a/V^{a+1})(M,F)$ and $(V^b/V^{b+1})(M,F)$ are viewed as
filtered left $\Dc_{Y'}$-module.
(For the shift of filtration $[-1]$, see (1.2.3).)

As for the right-hand side of (3.3.1), we have by definition
$$\psi_t(M,F)=\mopl_{a\in(0,1]}\Gr_V^a(M,F)[1],
\leqno(3.3.6)$$
where the shift of complex by 1 comes from Convention~4.
So we get (3.3.1) by using (1.2.2) and (3.3.4).
This finishes the proof of Proposition~(3.3).

\ms\nin
{\bf 3.4.~Proposition.} {\it Let $\M'(s'_1,\dots,s'_r)$ and
$\Lambda_y(T^*_{\vir}X)$ be as in the introduction and $(1.3)$
respectively. Then
$$\DR_y[\M'(s'_1,\dots,s'_r)]=(-1)^{\dim X}\Lambda_y(T^*_{\vir}X)
\q\h{in}\,\,\,K_0(X)[[y]][y^{-1}].$$
Hence $\Lambda_y(T^*_{\vir}X)\in
K_0(X)[y]\,\bl(=K_0(X)[y,y^{-1}]\cap K_0(X)[[y]]\br)$, and}
$$T_{y*}\bl(\M'(s'_1,\dots,s'_r)\br)=
(-1)^{\dim X}T^{\,\vir}_{y*}(X).$$

\ms\nin
{\it Proof.} It is enough to show the first equality, since the last
assertion follows from it by [Ve1], Th.~7.1 applied to the inclusion
$X\into\Zc_T$.
In the notation of the introduction, set
$$T_j:=\{t_k=0\,\,(k\le j)\}\subset T,\q\Zc_j:=\Zc_T\times_TT_j,$$
with the inclusions $i_j:\Zc_j\into\Zc_{j-1}$. Let $i_X:X\into\Zc_T$
denote the composition of the $i_j$ for $j\in[1,r]$. Set $n:=\dim X$.
Applying Proposition~(3.3) to the $i_j$, we get
$$(-1)^n\DR_y[\M'(s'_1,\dots,s'_r)]=
i^*_X\DR_y[\Q_{h,\Zc_T}](1+y)^{-r}.
\leqno(3.4.1)$$

Let $U$ be a non-empty Zariski-open subset of $T$ such that
$\Zc_U:=\Zc_T\times_TU$ is smooth over $U$. Let $T^*(\Zc_U/U)$ denote
the relative cotangent bundle. Since $T^*U$ is trivial, we have
$$\aligned &\DR_y[\Q_{h,\Zc_U}](1+y)^{-r}=
\Lambda_y[T^*\Zc_U](1+y)^{-r}\\
&=\Lambda_y\bl[T^*(\Zc_U/U)\br]=
i_{\Zc_U}^*\scc pr_1^*\bl(\Lambda_y[T^*Y^{(0)}]/
\mprod_{j=1}^r\Lambda_y[L_j^*]\br),\endaligned
\leqno(3.4.2)$$
where $i_{\Zc_U}:\Zc_U\into Y^{(0)}\times U$
and $pr_1:Y^{(0)}\times U\to Y^{(0)}$ are natural morphisms, and
$L_j^*$ is the dual of a very ample line bundle $L_j$ such that
$s_j\in\Ga(Y^{(0)},L_j)$.
(Here it is not necessary to assume that $L_j=L^{\otimes a_j}$.)

In order to apply inductively (2.4.3) and (2.4.4), we have to take
Zariski-open subsets $U'_j$, $U_j$ of $T_j\,\,(j\in[0,r-1])$
satisfying the two conditions
$$U'_j\setminus T_{j+1}=U_j\setminus T_{j+1},\q
U'_{j+1}:=U_j\cap T_{j+1}\ne\emptyset,
\leqno(3.4.3)$$
by increasing induction on $j\in[0,r-1]$, where $U'_r=T_r=\{0\}$, and
(2.4.3), (2.4.4) will be applied to the base changes of the inclusions
$$U'_{j+1}=U_j\cap T_{j+1}\into U_j.$$
If $j=0$, we set $U'_0:=U$. In general, if $U'_j$ is given, then
$U_j$ can be any Zariski-open subset of $T_j$ satisfying the two
conditions in (3.4.3). Here a canonical choice would be the maximal one
satisfying the two conditions, i.e.
$$U_j:=T_j\setminus\overline{(T_j\setminus(U'_j\cup T_{j+1}))}.$$
This means that $T_j\setminus U_j$ is the union of the irreducible
components of $T_j\setminus U'_j$ which are not contained in $T_{j+1}$.

Consider now the inclusions
$$\Zc_{U'_{j+1}}\into \Zc_{U_j},
\leqno(3.4.4)$$
which are obtained by the base change of the inclusions
$U'_{j+1}\into U_j\into U$ by $\Zc_U\to U$ for $j\in[0,r-1]$.
The assertion then follows from the `well-definedness' of (2.4.3) and
(2.4.4) (as is explained in (2.4)) which we apply to the inclusions
(3.4.4) inductively. This finishes the proof of Proposition~(3.4).

\newpage%\bs\bs
\centerline{\bf 4. Proof of the main theorem}
\bs\nin
In this section we prove Theorem~1 after showing Proposition~(4.1) below.
A proof of the latter has been presented in the proof of [PP], Prop.~5.1
in case $Y$ is smooth, where they demonstrated the contractibility of
the Milnor fiber using the integration of a controlled vector field as
in [Mi]. We give a sheaf-theoretic proof using an embedded resolution.
\ms\nin
{\bf 4.1.~Proposition.} {\it Let $X$ be a hypersurface of a
projective variety $Y$ such that $X$ is a very ample divisor, and
$Y\setminus X$ is smooth. Let $X'$ be a sufficiently general member
of the linear system associated with $X$ so that we have a
one-parameter family $\Y=\coprod_{c\in\PP^1}Y_c$ with $Y_0=X$,
$Y_{\infty}=X'$. Here $\PP^1$ is identified with $\C\cup\{\infty\}$ by
the inhomogeneous coordinate $t$ of $\PP^1$.
Let $\pi:\Y\to\PP^1$ be the projection. Set $X'':=X\cap X'$. Then}
$$\varphi_{\pi^*t}\Q_{\Y}|_{X''}=0,\q\h{\it i.e.}\q
\psi_{\pi^*t}\Q_{\Y}|_{X''}=\Q_{X''}.$$
\ms\nin
{\it Proof.}
We first treat the following simple case.
\ms\nin
{\bf (a) Normal crossing case.}
Assume that $Y$, $X'$ are smooth, $X$ is a divisor with simple normal
crossings, and $X'$ transversally intersects any intersections of
irreducible components of $X$.
(Here it is not necessary to assume that $X$ is a very ample divisor.)
Let $g$ be a local equation of $X$ around $x\in X''=X\cap X'$.
Let $y_1,\dots,y_n$ be local coordinates of $Y$ such that locally
$$g=\mprod_{i=1}^r\,y_i^{m_i},\q X'=\{y_n=0\},$$
where $r<n:=\dim Y$ and $m_i\ge 1\,(i\in[1,r])$.
Then we have locally
$$\Y=\{g=y_nt\}\subset Y\mtim\C,$$
where $t$ is identified with the affine coordinate of $\C\subset\PP^1$.
So the Milnor fiber of $\pi:\Y\to\PP^1$ around
$x\in X\cap X'\subset\pi^{-1}(0)$ is given by
$$\bl\{(y_1,\dots,y_n)\in\C^n\,\big|\,g=y_nt,\,\msum_{i=1}^n\,
|y_i|^2<\e^2\br\}.$$
where $0<|t|\ll\e\ll 1$.
This is identified with
$$\bl\{(y_1,\dots,y_{n-1})\in\C^{n-1}\,\big|\,\,\msum_{i=1}^{n-1}\,
|y_i|^2+|t|^{-2}\mprod_{i=1}^r\,|y_i|^{2m_i}<\e^2\br\}.$$
This set is contractible by using the natural action of
$\lambda\in[0,1]$ defined by
$$\lambda(y_1,\dots,y_{n-1})=(\lambda y_1,\dots,\lambda y_{n-1}).$$
So the assertion in the normal crossing case follows.
\ms\nin
{\bf (b) General case.}
Let $\si:(\Yt,\X)\to(Y,X)$ be an embedded resolution of singularities
such that $\X:=\si^{-1}(X)$ is a divisor with simple normal crossings.
By assumption, $\X':=\si^{-1}(X')$ is smooth and transversally
intersects any intersections of irreducible components of $\X$. So $\X'$
is the total transform of $X'$, and we get the one-parameter family
$$\Yct:=\mcoprod_{c\in\PP^1}\,\Yt_c\q\h{with}\q \Yt_0=\X,\,\,
\Yt_{\infty}=\X'.$$
Set
$$\Y_U:=\mcoprod_{c\in U}\,Y_c,\q\Yct_U:=\mcoprod_{c\in U}\,\Yt_c\q
\h{with}\q U:=\PP^1\setminus\{0\}.$$
Let $\pi_U:\Y_U\to\PP^{1}$, $\pit_U:\Yct_U\to\PP^{1}$ denote the canonical
morphisms. Then $\si$ induces
$$\sit_U:\Yct_U\to\Y_U\q\h{with}\q\pit_U=\pi_U\scc\sit_U:\Yct_U\to\Y_U\to
U.$$
Let $\si_c:\Yt_c\to Y_c\,(c\in U)$ be morphisms induced by $\si$. Set
$$\X'':=\X\cap\X',\q X'':=X\cap X'\q\h{with}\q\si'':=
\si|_{\X''}:\X''\to X''.$$
By the definition of $\Yct_U$, $\Y_U$, we have the inclusions
$$\iti''_U:\X''\times U\into\Yct_U,\q
i''_U:X''\times U\into\Y_U\q\h{over}\,\,\,U,$$
and $\sit$ induces an isomorphism
$$\Yct_U\setminus\iti''_U\bl(\X''\times U\br)\simto
\Y_U\setminus i''_U\bl(X''\times U\br).
\leqno(4.1.1)$$

Consider now the distinguished triangle
$$\Q_{\Y_U}\buildrel{\sit_U^{\#}}\over\too\R(\sit_U)_*\Q_{\Yct_U}
\to{\rm Cone}\,\sit_U^{\#}\buildrel{+1}\over\too.
\leqno(4.1.2)$$
By the definition of $\sit_U$ and using (4.1.1), we have
$${\rm Cone}\,\sit_U^{\#}\cong(i''_U)_*pr_1^*K''\q\h{with}\q
K'':={\rm Cone}\bl(\Q_{X''}\to\R\si''_*\Q_{\X''}\br),
\leqno(4.1.3)$$
where $pr_1:X''\times U\to X''$ is the first projection.

Let $t$ be the affine coordinate of $\C\subset\PP^1$.
Since the nearby cycle functor commutes with the direct image by
a proper morphism, we have a canonical isomorphism
$$\R(\si_0)_*\psi_{\pit^*t}\Q_{\Yct_U}=
\psi_{\pi^*t}\R\sit_*\Q_{\Yct_U}.$$
By the assertion in the normal crossing case together with the proper
base change theorem, we then get
$$\aligned &\R\si''_*\Q_{\X''}\simto
\R\si''_*\bl(\psi_{\pit^*t}\Q_{\Yct_U}\big|_{\X''}\br)\\
=\,\,&\R(\si_0)_*\psi_{\pit^*t}\Q_{\Yct_U}\big|_{X''}
=\psi_{\pi^*t}\R\sit_*\Q_{\Yct_U}\big|_{X''}.\endaligned
\leqno(4.1.4)$$
Let $i_0:Y_0\into\Y$ denote the inclusion. Apply the functorial
morphism $i_0^*\to\psi_{\pi^*t}$ to the distinguished triangle
$$\Q_{\Y}\buildrel{\sit^{\#}}\over\too\R(\sit)_*\Q_{\Yct}
\to{\rm Cone}\,\sit^{\#}\buildrel{+1}\over\too,$$
which is the extension of (4.1.2) over $\PP^1$.
Restricting these over $X''\subset X=Y_0$, and using the
`well-definedness' of the nearby cycle functor as in (2.4.4), we then
get a morphism of distinguished triangles
$$\begin{matrix}\Q_{X''}&\to&\R\si''_*\Q_{\X''}&\to&K''&
\buildrel{+1}\over\too\\
\downarrow\,&&\downarrow&&||\,\\
\psi_{\pi^*t}\Q_{\Y_U}|_{X''}&\to&
\psi_{\pi^*t}\R\sit_*\Q_{\Yct_U}|_{X''}&\to&K''&
\buildrel{+1}\over\too\end{matrix}
\leqno(4.1.5)$$
where the right vertical isomorphism follows from (4.1.3).
Moreover, the middle vertical morphism is an isomorphism by (4.1.4).
So the left vertical morphism is an isomorphism.
This finishes the proof of Proposition~(4.1).

\ms\nin
{\bf 4.2.~Proof of Theorem~1.}
We proceed by increasing induction on $r$.
We fix $s_1,\dots,s_r$ such that
$$\Si:=\Sing\,X\supset\Sing\,Y\supset\Si':=\Sing\,X'=\Sing\,Y\cap X',$$
in the notation of the introduction, e.g.
$Y:=Y^{(r-1)}$, $X:=Y^{(r)}$, $X'=Y\cap s'_r{}^{-1}(0)$.

Let $\Y$ be the blow-up of $Y$ along $X'':=X\cap X'$ so that
$$\Y=\mcoprod_{c\in\PP^1}\,Y_c\q\h{with}\q Y_0=X,\,\,Y_{\infty}=X'.$$
Here we use $c$ to denote a point of $\PP^1$, and the coordinate $t_r$
will be used to define the nearby and vanishing cycle functors.

Set $T_r:=\{t_i=0\,\,(1\le i<r)\}\subset T$. Then
$$\Y\setminus Y_{\infty}\cong\Zc_T\cap(Y^{(0)}\times T_r).$$
Set $\Y_U:=\Y\times_{T_r}U$ with $U$ a sufficiently small non-empty
Zariski-open subset of $T_r$ such that the hypersurface
$\{s_r=c\,s'_r\}\subset Y^{(0)}$ is smooth and intersects
transversally any strata of Whitney stratifications of
$Y^{(k)}\,\,(0<k<r)$ for any $c\in U$.
Set $U':=U\cup\{\infty\}\subset\PP^1$.
Let $i_c:Y_c\into\Y_{U'}$ denote the inclusion for $c\in U'$.
Set
$$\M'_{U'}:=\M'(s'_1,\dots,s'_{r-1})|_{\Y_{U'}},\q
\M'_c:=i_c^*\M'_{U'}[-1]\,\,(c\in {U'}).$$
These are mixed Hodge modules on $\Y_{U'}$ and $Y_c\,\,(c\in U')$
respectively, and moreover $Y_c$ transversally intersects any
strata of a Whitney stratification of $\Y_U$ compatible with $\M'_U$
by shrinking $U$ if necessary.
(Here $U'$ may contain $\infty\in\PP^1$, since $X'=\{s'_r=0\}\subset Y$
is assumed to be sufficiently general in the one-parameter family
$Y_c=\{s_r=c\,s'_r\}\subset Y\,\, (c\in U)$.)

Set $n=\dim X$.
By inductive hypothesis applied to the smooth projective variety
$\{s_r=c\,s'_r\}\subset Y^{(0)}$ for $c\in U'$, we have short exact
sequences of mixed Hodge modules on $Y_c\,\,(c\in U')$
$$0\to\Q_{h,Y_c}[n]\to\M'_c\to\M_c\to 0,
\leqno(4.2.1)$$
where the last term is defined by the cokernel of the injection.
(Indeed, the restriction by the inclusion
$\{s_r=c\,s'_r\}\into Y^{(0)}$ commutes with the nearby cycle functors
by [DMST].)
They imply a short exact sequence of mixed Hodge modules on $\Y_{U'}$
$$0\to\Q_{h,\Y_{U'}}[n+1]\to\M'_{U'}\to\M_{U'}\to 0,
\leqno(4.2.2)$$
together with the isomorphism
$$\M_c=i_c^*\M_{U'}[-1]\,\,(c\in U').$$
Applying the nearby cycle functor $\psi_{t_r}$ to (4.2.2), we get an
exact sequence of mixed Hodge modules on $X$
$$0\to\psi_{t_r}\Q_{h,\Y_{U'}}[n]\to\M'(s'_1,\dots,s'_r)\to
\psi_{t_r}\M_{U'}[-1]\to 0.
\leqno(4.2.3)$$
Here we use the well-definedness of (2.4.4) which implies that
$\psi_{t_r}$ can be defined for mixed Hodge modules defined on the
complement of $\{t_r=0\}$.
We have moreover the short exact sequence
$$0\to\Q_{h,X}[n]\to\psi_{t_r}\Q_{h,\Y_{U'}}[n]\to
\varphi_{t_r}\Q_{h,\Y}[n]\to 0,
\leqno(4.2.4)$$
since $X$ is a complete intersection so that $\Q_X[n]$ is a
perverse sheaf.
These exact sequences imply the injectivity of the natural morphism
$$\Q_{h,X}[n]\to\M'(s'_1,\dots,s'_r).$$
Let $\M(s'_1,\dots,s'_r)$ be its cokernel. Then we have a short exact
sequence of mixed Hodge modules
$$0\to\varphi_{t_r}\Q_{h,\Y}[n]\to\M(s'_1,\dots,s'_r)\to
\psi_{t_r}\M_{U'}[-1]\to 0.
\leqno(4.2.5)$$

Set $f:=(s_r/s'_r)|_{Y\setminus X'}$ as in the introduction.
Since $\Si:=\Sing\,X\supset\Sing\,Y$, the support of
$\varphi_f\Q_{h,Y}$ is contained in $\Si\setminus X'$, and
$\varphi_f\Q_{h,Y}[n]$ can be viewed as a mixed Hodge module on
$\Si\setminus X'$. By Proposition~(4.1) we get
$$\varphi_{t_r}\Q_{h,\Y}[n]=(i_{\Si\setminus X',\Si})_{!\,}
\varphi_f\Q_{h,Y}[n],
\leqno(4.2.6)$$
where $i_{X\setminus X',X}:X\setminus X'\into X$ is the inclusion,
and the left-hand side is viewed as a mixed Hodge module on $\Si$.
So the proof of Theorem~1 is reduced to Proposition~(4.3) below,
except for the independence of the $s'_j$.
To show the latter, we consider the one-parameter families
$$s'_{j,\lambda}:=(1-\lambda)s'_j+\lambda s''_j\,\,\,(\lambda\in\C)
\q\h{for}\,\,\,j\in[1,r],$$
if we are given sufficiently general $s'_j$ and $s''_j$ for $j\in[1,r]$,
where $\lambda$ is independent of $j\in[1,r]$.
Then the assertion follows from Lemma~(2.6).
(See also the proof of Proposition~(4.4) below.)
\ms\nin
{\bf 4.3.~Proposition.} {\it With the above notation,
$\psi_{t_r}\M_{U'}[-1]$ is supported in $\Si'$, and}
$$M_y(X')=(-1)^nT_{y*}\bl(\psi_{t_r}\M_{U'}[-1]\br).$$

\ms\nin
{\it Proof.} By inductive hypothesis, we have
$${\rm supp}\,\M_c\subset\Sing\,Y_c=\Si'\q(c\in U'),$$
where the last equality holds since $s'_r$ is assumed to be sufficiently
general.
So $\M_{U'}$ is viewed as a mixed Hodge module on $\Si\times U'$, and
this can be extended to a mixed Hodge module on $\Si\times\C$ by
extendability of mixed Hodge modules in the algebraic case, where
$\C$ is the complement of some point of $\PP^1$.
By Proposition~(3.3) together with (2.4.2) and Lemma~(2.6),
we then get the assertion.
(Note that the non-characteristic restriction $i^*_c$ for any $c\in U'$
trivially commutes with $\DR_y$ up to the factor as in Proposition~(3.3).)
This finishes the proofs of Proposition~(4.3) and Theorem~1.

\ms\nin
{\bf 4.4.~Proposition.} {\it In Theorem~$1$, $M_y(X)$ is independent
of the choice of the $s_j$ if $([s_j])$ belongs to a sufficiently
small non-empty Zariski-open subset $\U_I$ of
$\prod_{j=1}^r\PP(I_{X,a_j})$.}

\ms\nin
{\it Proof.} Since $\M(s'_1,\dots,s'_r)$ in the introduction may depend
also on $s_1,\dots,s_r$, we denote it by $\M(\sss,\sss')$ with
$\sss:=(s_1,\dots,s_r)$, $\sss':=(s'_1,\dots,s'_r)$. By the construction
of $\M(\sss,\sss')$ in the introduction, there is a non-empty Zariski-open
conical subset $\U_{I,R}$ of
$$\mprod_{j=1}^rI_{X,a_j}\times\mprod_{j=1}^rR_{L,a_j},$$
together with a mixed Hodge module $\M_{\U_{I,R}}$ on $\Si\times\U_{I,R}$
such that $\M(\sss,\sss')$ is the restriction of $\M_{\U_{I,R}}$ by the
inclusion $\Si\times\{\sss,\sss'\}\into\Si\times\U_{I,R}$, and moreover
this inclusion is strictly non-characteristic for the underlying perverse
sheaf $\F_{\U_{I,R}}$ of $\M_{\U_{I,R}}$, i.e.\ $\F_{\U_{I,R}}$ is a
topologically locally constant family of perverse sheaves parametrized
by $\U_{I,R}$. (This can be shown for instance by using [BMM] together
with a Whitney stratification of the total space, since its restriction
to a sufficiently general fiber is a Whitney stratification.)

We then get the independence of $M_y(X)$ by the choices of
$\sss,\sss'$ using Lemma~(2.6) together with a one-parameter family
as in the last part of the proof of Theorem~1 in (4.2).
Since the independence of $\sss'$ is already shown, the assertion follows.

\ms\nin
{\bf 4.5.~Remarks.} (i) Let $\Zc_T\to T$ be as in the introduction.
If $X=\Zc_0$ has an isolated singular point $x$, then we have a mixed
Hodge structure on the vanishing cohomology. This can be defined by
restricting the morphism $\Zc_T\to T$ over a generic smooth curve on $T$ passing through $0$ and
using an embedded resolution of singularities as in [St2]. (This can be
calculated also by using mixed Hodge modules as in [BS], Th.~4.3.)
The associated mixed Hodge numbers are independent of the choice of the
generic smooth curve by the theory of spectra for arbitrary varieties
which uses the deformation to the normal cone (see e.g.\ [DMS],
Remark~(1.3)(i)).
By a similar argument these mixed Hodge numbers coincide with those
given by Corollary~2 assuming that $s_1,\dots,s_{r-1}$ and
$s'_1,\dots,s'_r$ are sufficiently general (by using Proposition~(3.2)).

\ms
(ii) For a germ of a complete intersection with an isolated singularity
$(X,0)$, there is a versal flat deformation $\rho:(\Xc,0)\to(S,0)$ in
the following sense: any flat deformation $(\Xc',0)\to(S',0)$ of $(X,0)$
is analytically isomorphic to the pull-back of $\rho$ by some complex
analytic morphism $(S',0)\to(S,0)$, see [KS], [Tju].
Choosing a generic smooth curve on $S$ passing through $0$, we also get
a mixed Hodge structure on the vanishing cohomology, and the associated
mixed Hodge numbers are independent of the choice of the generic curve,
see also [ES]. Moreover these numbers coincide with those given in
Remark~(i) above by using [DMS], Cor.~3.4.

\bs\bs\centerline{\bf 5. Explicit calculations of the Hirzebruch-Milnor
classes}
\bs\nin
In this section we show how the Hirzebruch-Milnor classes can be
calculated explicitly.
\ms\nin
{\bf 5.1.~Stratification.} 
The first term of the right-hand side of (0.2) in Theorem~1 can be
described explicitly as follows.
Let $\Sc$ be a complex algebraic stratification of $\Si\setminus X'$
such that the $\Hc^j\varphi_f\Q_Y|_S$ are local systems for any
$j\in\N$ and $S\in\Sc$.
Then $\Hc^j\varphi_f\Q_Y|_S$ underlies the variation of mixed
Hodge structure $\Hc^ji_{S,\Si\setminus X'}^*\varphi_f\Q_{h,Y}$.
Here it is enough to assume that each stratum $S\in\Sc$ is a
locally closed smooth subvariety and the Whitney conditions are
not needed as long as we have local systems on each stratum.
We have the following.
\ms\nin
{\bf 5.2.~Proposition.} {\it With the above notation and assumption,
we have}
$$T_{y*}\bl((i_{\Si\setminus X',\Si})_{!\,}\varphi_f\Q_{h,Y}\br)=
\msum_{j\in\N,S\in\Sc}\,(-1)^j\,T_{y*}\bl((i_{S,\Si})_!
(\Hc^ji_{S,\Si\setminus X'}^*\varphi_f\Q_{h,Y})\br).
\leqno(5.2.1)$$
\ms
Indeed, this follows from the following.
\ms\nin
{\bf 5.3.~Proposition.} {\it Let $X$ be a complex algebraic variety.
Let $\MM\in D^b\MHM(X)$, and $K^{\ssb}$ its underlying $\Q$-complex.
Let $\Sc=\{S\}$ be a complex algebraic stratification of $X$ such
that for any $S\in\Sc$, $S$ is smooth, $\So\setminus S$ is a union of
strata, and the $\Hc^iK^{\ssb}|_S$ are local systems on $S$ for any $i$.
Let $j_S:S\into X$ denote the inclusion. Then
$$T_{y*}(\MM)=\msum_{S,i}\,(-1)^i\,T_{y*}\bl((j_S)_!\Hc^i
(j_S)^*\MM\br),
\leqno(5.3.1)$$
where $\Hc^i$ for $(j_S)^*\MM$ is associated to the classical
$t$-structure as in {\rm [Sa2], 4.6.2} $($which coincides with the usual
$t$-structure on the derived category of mixed Hodge modules up to a
shift if it is restricted to a stratum of the stratification$)$ so that
$\Hc^iK^{\ssb}|_S$ underlies variations of mixed Hodge structures
$\Hc^i(j_S)^*\MM$.}
\ms\nin
{\it Proof.} Set $X_k:=\bigcup_{\dim S\le k}S\subset X$.
We have the canonical inclusions
$$j_k:X_k\setminus X_{k-1}\into X_k,\q i_k:X_k\into X,$$
together with the distinguished triangles
$$(i_k)_*(j_k)_!j_k^*i_k^*\MM\to(i_k)_*i_k^*\MM\to
(i_{k-1})_*i_{k-1}^*\MM\buildrel{+1}\over\to.
\leqno(5.3.2)$$
These imply the following identity in the Grothendieck group of
$D^b\MHM(X)$
$$[\MM]=\msum_k\,\bl[(i_k)_*(j_k)_!j_k^*i_k^*\MM\br].
\leqno(5.3.3)$$
Moreover, $j_S$ is the composition of $j'_S:S\into X_k$
with $i_k$ where $k=\dim S$.
So the assertion follows.
\ms\nin
{\bf 5.4.~Logarithmic forms.} 
The right-hand side of (5.2.1) can be described as follows.
Let $\M$ be an admissible variation of mixed Hodge structure on a
stratum $S$ more generally, where we may assume that $\M$ is a
polarizable variation of Hodge structure by taking the graded pieces
of the weight filtration $W$.
Let $M$ be the underlying $\Oc_S$-module with the Hodge filtration
$F$ of the variation of mixed Hodge structure.
Take a smooth partial compactification $i_{S,Z}:S\into Z$ such that
$D:=Z\setminus S$ is a divisor with simple normal crossings and
moreover $i_{S,\Si}=\pi_Z\scc i_{S,Z}$ for a proper morphism
$\pi_Z:Z\to\Si$.
Here $D$ cannot be empty, since we take a stratification
of $\Si\setminus X'$ and $X'$ is a hyperplane section.
Let $M_Z^{>0}$ be the Deligne extension with eigenvalues of the
residues of the logarithmic connection contained in $(0,1]$.
Then we have the following.
\ms\nin
{\bf 5.5.~Proposition.} {\it With the above notation, we have}
$$T_{y*}\bl((i_{S,\Si})_!\M\br)=
\msum_{p,q}(-1)^q(\pi_Z)_*td_{(1+y)*}\bl[\Gr_F^pM_Z^{>0}
\otimes\Omega_Z^q(\log D)\br](-y)^{p+q}.
\leqno(5.5.1)$$
\ms\nin
{\it Proof.}
Let $(M_Z,F)$ be the underlying filtered $\Dc_Z$-module of
$(i_{S,Z})_!\M$.
Then we have a canonical inclusion
$$(M_Z^{>0}\otimes\Omega^{\ssb}_Z(\log D),F)[\dim Z]\into
\DR_Z(M_Z,F),
\leqno(5.5.2)$$
where the filtration $F$ on $(M_Z^{>0}\otimes\Omega^p_Z(\log D)$
is shifted by $-p$ as usual.
Moreover, it is a filtered quasi-isomorphism by [Sa2], Prop.~3.11,
see also [CMSS], 3.3.
So the assertion follows from the commutativity of $td_*$ with the 
pushforward by proper morphisms [BFM], see also [CMSS], 3.3.
\ms\nin
{\bf 5.6.~Motivic Milnor fibers.} 
It is also possible to use the theory of motivic nearby fibers [DL]
(see also [St1]) to calculate the first term of the right-hand side of
(0.2).
Let $\si:(\Yt,\X)\to(Y,X)$ be an embedded resolution inducing the
isomorphism outside the singular locus of $X_{\rm red}$ and such that
$\X:=\si^{-1}(X)$ and $\si^{-1}(\Si)$ are divisors with
simple normal crossings, where simple means that the irreducible
components $E_i$ of $\X$ are smooth.
(Such a condition for $\si^{-1}(\Si)$ is satisfied if
it is a divisor and the other conditions are satisfied.)
Let $m_i$ be the multiplicity of $\X$ along $E_i\,(i=1,\dots,s)$.
Put $\X':=\si^{-1}(X')$.
For $I\subset\{1,\dots,s\}$, set
$$E_I:=\mcap_{i\in I}E_i,\q E'_I:=E_I\setminus\X',\q
E^{\pc}_I:=E'_I\setminus\bl(\mcup_{i\notin I}E_i\br).$$
Let $\f$ be the pull-back of $f$ to $\Yt\setminus\X'$.
There are smooth varieties $\E^{\pc}_I$ together with finite
\'etale morphisms $\gamma_I:\E^{\pc}_I\to E^{\pc}_I$ such that
$(\gamma_I)_*\Q_{\E^{\pc}_I}=\Hc^0\psi_{\f}\Q_{\Yt}$, see [DL].
Note that we can get these varieties by taking the normalization of
the base change of $\f$ by a ramified $m$-fold covering of an open
disk where $m={\rm LCM}(m_i)$, see [St1].
They may be called the Stein factorization of the Milnor fibers,
see [Lo].
Let $\jti_I:\E^{\pc}_I\into\E_I$ be a smooth compactification such that
$\gamma_I$ is extended to $\g_I:\E_I\to E_I$.
We have canonical morphisms
$$i_{I,J}:\E_{I,J}\into\E_I,\q\si_{\Si}:=\si|_{\si^{-1}(\Si)}:
\si^{-1}(\Si)\to X.$$
Set $\La':=\{I\subset\{1,\dots,r\}\mid E_I\subset\si^{-1}(\Si)\}$.
Then we have the following.
\ms\nin
{\bf 5.7.~Proposition.} {\it With the above notation and assumption,
we have}
$$\aligned T_{y*}\bl((i_{\Si\setminus X',\Si})_{!\,}\varphi_f\Q_{h,Y}
\br)&=\msum_{I\in\La'}\,(\si_{\Si}\scc\g_I)_*T_{y*}
\bl((\jti_I)_{!\,}\Q_{h,\E^{\pc}_I})\br)(1+y)^{|I|-1}\\
&\q-T_{y*}(\Si)+(i_{\Si\cap X',\Si})_*T_{y*}(\Si\cap X').\endaligned
\leqno(5.7.1)$$
\ms\nin
{\it Proof.}
Let $h:\Zc'\to\C$ be the normalization of the base change of
$\f:\Yt':=\Yt\setminus\X'\to\C$ by the $m$-fold ramified covering
$\C\to\C$ defined by $t\mapsto t^m$ where $m={\rm LCM}(m_i)$ with
$m_i$ the multiplicity of $\X$ along $E_i$. Set $\Zc'_0:=h^{-1}(0)$.
There is a canonical morphism $\gamma:\Zc'_0\to\X\setminus\X'$. Set
$$\E_I^{\pc}:=\gamma^{-1}(E_I^{\pc}),\q
m_I:={\rm GCD}\bl(m_i\,(i\in I)\br).$$
Then the induced morphism $\gamma_I:\E_I^{\pc}\to E_I^{\pc}$
is finite \'etale of degree $m_I$.
By [St1] we have
$$\Hc^i\psi_h\Q_{\Zc'}|_{\E_I^{\pc}}=\buildrel{\nu(I,i)}\over\mopl
\Q_{\E_I^{\pc}}(-i)\q\h{with $\,\,\nu(I,i)=\binom{|I|-1}{i}$},$$
since $\Zc'_0$ is a divisor with $V$-normal crossings on a
$V$-manifold $\Zc'$. So we get
$$\Hc^i\psi_{\f}\Q_{\Yt'}|_{E_I^{\pc}}=\buildrel{\nu(I,i)}\over\mopl
(\gamma_I)_*\Q_{\E_I^{\pc}}(-i).$$
This is an isomorphism as $\Q$-complexes.
However, it holds as variations of Hodge structures of type $(i,i)$
since the Hodge filtration is trivial.

Let $j:X\setminus X'\into X$ denote the inclusion.
By Proposition~(4.1), we have
$$\varphi_{\pi^*t}\Q_{h,\Y}=j_!\varphi_f\Q_{h,Y\setminus X'}=
(i_{\Si})_*\,i_{\Si}^*\,j_!\bl(C(\Q_{h,X\setminus X'}\to
\psi_f\Q_{h,Y\setminus X'})\br).$$
Here the second isomorphism follows from the distinguished triangle
$$\Q_{h,X\setminus X'}\to\psi_f\Q_{h,Y\setminus X'}\to
\varphi_f\Q_{h,Y\setminus X'}\buildrel{+1}\over\longrightarrow,$$
together with the fact that $\varphi_f\Q_{h,Y\setminus X'}$ is supported
on $\Si\setminus X'$.

By the commutativity of the nearby cycle functor with the direct
image by a proper morphism, we have
$$\psi_f\Q_{h,Y\setminus X'}=
(\si'_0)_*\psi_{\f}\Q_{h,\Yt\setminus\X'},$$
where $\si'_0:\X\setminus\X'\to X\setminus X'$ is the restriction
of $\si$.
So the assertion follows from Proposition~(5.3).
This finishes the proof of Proposition~(5.7).
\ms\nin
{\bf 5.8.~Remarks.} (i)
The term $(1+y)^{|I|-1}$ in (5.7.1) comes from the mixed Hodge structure
on the cohomology of the irreducible components of the Milnor fiber
which is homeomorphic to $(\C^*)^{|I|-1}$ in the normal crossing case
and corresponds to $(1-\LL)^{|I|-1}$ in [DL].
(Note that the stalks of the nearby cycle sheaves are given by the
cohomology of the Milnor fibers.)
Here we may assume that $\E_I\setminus\E^{\pc}_I$ is a divisor with
simple normal crossings replacing $\E_I$ if necessary.
In this case, let $\E_{I,j}$ be the irreducible components of
$\E_I\setminus\E^{\pc}_I$ which are assumed smooth.
Set $\E_{I,J}:=\mcap_{j\in J}\E_{I,j}$ where $\E_{I,\emptyset}=\E_I$
if $J=\emptyset$. Then the usual resolution argument implies
$$T_{y*}\bl((\jti_I)_{!\,}\Q_{h,\E^{\pc}_I}\br)=
\msum_J\,(-1)^{|J|}(i_{I,J})_*T_{y*}(\E_{I,J}).
\leqno(5.8.1)$$

(ii) In Theorem~1 and its corollaries and also in Propositions
of this section, most of the formulas hold also on the level of the
Grothendieck group of mixed Hodge modules and hence on that of coherent
sheaves by using $\DR_y$ (except for (5.5.1) which is also valid, or
rather meaningful, only on the level of the Grothendieck group of
coherent sheaves).

\end{document}